\documentclass[10pt,a4paper]{amsart}

\usepackage{amssymb}

\input xy
\xyoption{all}
\usepackage{tikz}
\usetikzlibrary{arrows,shapes}

\renewcommand{\epsilon}{\varepsilon}
\renewcommand{\setminus}{\smallsetminus}
\renewcommand{\emptyset}{\varnothing}

\newtheorem{theorem}{Theorem}[section]
\newtheorem{proposition}[theorem]{Proposition}
\newtheorem{corollary}[theorem]{Corollary}
\newtheorem{lemma}[theorem]{Lemma}
\newtheorem{conjecture}[theorem]{Conjecture}

\theoremstyle{definition}
\newtheorem{example}[theorem]{Example}

\newtheorem{definition}[theorem]{Definition}

\theoremstyle{remark}
\newtheorem{remark}[theorem]{Remark}

\newcommand{\XX}{\mathfrak X}

\newcommand{\FF}{\mathfrak F}
\newcommand{\FK}{{{\mathfrak F}_k}}

\newcommand{\Q}{\mathbb Q}
\newcommand{\Z}{\mathbb Z}
\newcommand{\N}{\mathbb N}

\newcommand{\fpinfty}{{\FP}_{\infty}}

\newcommand{\FFF}{\operatorname{F}}
\newcommand{\UF}{\underline{\operatorname{F}}}
\newcommand{\FP}{\operatorname{FP}}

\newcommand{\UFP}{\underline{\operatorname{FP}}}

\newcommand{\EXT}{\operatorname{Ext}}
\newcommand{\TOR}{\operatorname{Tor}}

\newcommand{\Ho}{\operatorname{H}}
\newcommand{\cohom}[3]{H^{{\raise1pt\hbox{$\scriptstyle#1$}}}(#2\>\!,#3)}
\newcommand{\tatecohom}[3]%
  {\widehat H^{{\raise1pt\hbox{$\scriptstyle#1$}}}(#2\>\!,#3)}

\newcommand{\Cohom}[3]%
  {H^{{\raise1pt\hbox{$\scriptstyle#1$}}}\big(#2\>\!,#3\big)}
\newcommand{\Tatecohom}[3]%
  {\widehat H^{{\raise1pt\hbox{$\scriptstyle#1$}}}\big(#2\>\!,#3\big)}

\newcommand{\homol}[3]{H_{{\lower1pt\hbox{$\scriptstyle#1$}}}(#2\>\!,#3)}
\newcommand{\homolog}[2]{H_{{\lower1pt\hbox{$\scriptstyle#1$}}}(#2)}

\newcommand{\colim}{\varinjlim}
\newcommand{\invlim}{\underset{\longleftarrow}{\lim}}

\newcommand{\mono}{\rightarrowtail}
\newcommand{\epi}{\twoheadrightarrow}

\renewcommand{\implies}{\Rightarrow}

\newcommand{\eg}{{\underline EG}}

\newcommand{\EXG}{E_\XX G}
\newcommand{\EFG}{E_\FF G}

\newcommand{\Hom}{\operatorname{Hom}}

\newcommand{\Tor}{\operatorname{Tor}}
\newcommand{\blah}{-}

\newcommand{\Mod}{\mathfrak{Mod}}

\newcommand{\OXG}{\mathcal O_{\mathcal X}G}
\newcommand{\OFG}{\mathcal O_{\mathcal F}G}

\newcommand{\OFKG}{\mathcal O_{{\mathcal F}_k}G}

\newcommand{\uz}{\underline\Z}

\title[Generalised Thompson groups]{Bredon cohomological finiteness conditions for generalisations of Thompson groups}\thanks{to appear Groups, Geometry and Dynamics}

\author{C.~Mart\'inez-P\'erez}

\address{Conchita Mart\'inez-P\'erez, Departamento de Matem\'aticas, Universidad de Zaragoza,
50009 Zaragoza, Spain} \email{conmar@unizar.es}

\author{B.~ E.~A.~Nucinkis}

\address{Brita E.~A.~Nucinkis,Department of Mathematics, Royal Holloway, University of London, Egham, TW20 0EX , United Kingdom}\email{brita.nucinkis@rhul.ac.uk}

\date{\today} 

\thanks{The first named author was partially supported by BFM2010-19938-C03-03, Gobierno de Arag\'on and European Union's ERDF funds.}

\begin{document}

\thispagestyle{empty}

\begin{abstract}  We define a family of groups that generalises Thompson's groups $T$ and $G,$ and also those of Higman, Stein and Brin. For groups in this family we describe  centralisers of finite subgroups and show, that for a given finite subgroup $Q$, there are finitely many  conjugacy classes  of finite subgroups isomorphic to $Q$. We consider  a slightly weaker property, quasi-$\underline\FFF_\infty$, to that of a group possessing a finite type model for the classifying space for proper actions $\eg,$ and give criteria for the T versions of our groups to be of type quasi-$\underline\FFF_\infty.$  We also generalise some well-known properties of ordinary cohomology to Bredon cohomology. 
\end{abstract}

\maketitle

\section{Introduction}

\medskip\noindent Thompson's groups $F$, $T$ and $G$ (also denoted $V$), which can be defined as certain homeomorphism groups of the unit interval, the circle and the Cantor-set, respectively, have received a large amount of attention in recent years. There are  many interesting generalisations  of these groups, such as the Higman-Thompson groups $F_{n,r}$, $T_{n,r}$, $G_{n,r}$ (recall that $T=T_{2,1}$ and $G=G_{2,1}$),  the $T$- and $G$-groups defined by Stein \cite{stein} and the higher dimensional Thompson groups $sV=sG_{2,1}$ defined by Brin \cite{Brin}. All these groups contain free abelian groups of infinite rank, are finitely presented and  with the exception of $sV$ for $s\geq 4,$ are known to be  of type $\FP_\infty$ \cite{brown, stein, desiconbrita2, bhm}.  Furthermore, the $G$- and $T$-groups contain finite groups of arbitrarily large orders. In this paper we  consider automorphism groups of certain Cantor algebras which include Higman-Thompson, Stein and Brin's groups.

As in the original exposition by Higman \cite{higman} and in Brown's proof \cite{brown} that $F_{n,r},$ $T_{n,r}$ and $G_{n,r}$ are of type $\FP_\infty,$ we consider a Cantor algebra $U_r(\Sigma)$ on a so called valid set of relations $\Sigma$ and define groups $G_r(\Sigma)$ as follows: the elements of $G_r(\Sigma)$ are  bijections between certain subsets of $U_r(\Sigma),$ which we call admissible. One can show that these groups are finitely generated, see \cite{mart11}. Provided that the relations in $\Sigma$ are order preserving we can also define the groups $T_r(\Sigma)$, which are  given by cyclic order preserving bijections. One can also define generalisations of $F_{n,r}$. 

 The admissible subsets of $U_r(\Sigma)$ form a poset, and the groups $T_r(\Sigma)$ and $G_r(\Sigma)$ act on the geometric realisation $|\frak A_r(\Sigma)|$ of this poset (for the original Thompson-Higman algebras this was already used by Brown in \cite{brown}).

\medskip\noindent  Let $G$ be either $T_r(\Sigma)$ or $G_r(\Sigma).$ For every finite subgroup $Q$ we consider the fixed point sets  $\frak A_r(\Sigma)^Q.$  The $Q$-set structure of every admissible subset  $Y\in \frak A_r(\Sigma)^Q$ is determined by its decomposition into transitive $Q$-sets.  We show (Theorem \ref{conjclasses})  that there are finitely many conjugacy classes in $G$ of subgroups isomorphic to $Q$. Furthermore we show (Theorem \ref{centralizer}) that there is an extension
$$K \mono C_{G_r(\Sigma)}(Q) \twoheadrightarrow G_{r_1}(\Sigma)\times\ldots\times G_{r_t}(\Sigma)$$
with locally finite kernel, where the $r_1,...,r_t$ are integers uniquely determined by $Q$. We also get the analogous result (Theorem \ref{Textension}) for the groups $T_r(\Sigma)$ (if defined) that, for a certain $l$ also determined by $Q,$ there is a central extension
$$Q \mono C_{T_r(\Sigma)}(Q) \epi T_{l}(\Sigma).$$
 This generalises a result of Matucci \cite[Theorem 7.1.5]{matuccithesis}  for the original Thompson group $T$.

\noindent  Recently a variant of the Eilenberg-Mac Lane space, the classifying space with respect to a family $\XX$ of subgroups, has been well researched.
Let $X$ be a $G$-CW-complex. $X$ is said to be a model for $\EXG,$ the classifying space with isotropy in the family $\XX$ if $X^K$ is contractible for $K \in\XX$ and $X^K$ is empty otherwise.  The classifying space $X$ for a family satisfies the following universal property: whenever there is a $G$-CW-complex $Y$ with isotropy lying in the family $\XX$, there is a $G$-map $Y \to X$, which is unique up to $G$-homotopy. In particular, $\EXG$ is unique up to $G$-homotopy equivalence.

\noindent
For the family $\mathfrak F$ of finite subgroups we denote $\EFG$ by $\eg$, the classifying space for proper actions. We say a group is of type $\UF_\infty$ if it admits a finite type model for $\eg.$ We show:

\medskip\noindent {\bf Theorem \ref{model}.} {\sl $|\frak A_r(\Sigma)|$ is a model for $\underline{E}G_r(\Sigma).$}

\medskip\noindent Obviously, this model has infinite dimension. Since these groups contain free abelian groups of infinite rank, they cannot possess any finite dimensional model.  Exactly as ordinary classifying spaces yield free resolutions which can be used to define ordinary group cohomology, classifying spaces with isotropy in a family produce free resolutions in a functor  category, which are used to define  Bredon cohomology. We shall review  properties of Bredon cohomology in Section 2. Many notions from ordinary cohomology have a Bredon analogue. For example, we say a group $G$ is of type Bredon-$\FP_\infty$ if there is a Bredon-projective resolution of the constant Bredon-module $\Z(-)$ by finitely generated Bredon-projective modules. The connection to classifying spaces and to ordinary cohomology is given by the following two results:

\begin{theorem}\cite[Theorem 0.1]{lm}\label{lueckmeintrupp}
A group $G$ has a finite type model for a classifying space with isotropy in a family if and only if
the group is of type Bredon-$\fpinfty$ and there is a model for the classifying space with finite $2$-skeleton.
\end{theorem}

\noindent In particular we say a group is of type $\UFP_\infty$ if it is of type Bredon-$\FP_\infty$ for the family of finite subgroups.

\begin{theorem}\cite[Theorem 4.2]{lueck}\label{lueckstheorem}
A group $G$ admits a finite type model for $\eg$ if and only $G$ has finitely many conjugacy classes of finite subgroups and for each finite subgroup $K$ of $G$ the centraliser $C_G(K)$ is of type $\fpinfty$ and finitely presented.
\end{theorem}

Equivalently, $G$ admits a finite type model for $\eg$ if and only it is of type $\UFP_\infty$ and centralisers of finite subgroups are finitely presented. The purpose of this paper is to study the possible finiteness conditions a model for $\eg$ for the groups $G_r(\Sigma)$ and $T_r(\Sigma)$ can satisfy.
Since the groups we are considering do not have a bound on the orders of their finite subgroups, we need to weaken the condition on the number of conjugacy classes. We consider the property  quasi-$\UFP_\infty,$ which has the same condition on the centralisers of finite subgroups as $\UFP_\infty$ but just requires that for each finite subgroup $Q$ of $G,$ there are only finitely many conjugacy classes of subgroups isomorphic to $Q$. Note, that for groups with a  bound on the orders of their finite subgroups both properties coincide. In \cite{learynucinkis} it was shown that there are examples of groups of type $\FP_\infty$, which have a bound on the orders of the finite subgroups, yet are not of type $\UFP_\infty.$ These examples are virtually torsion free, admit a finite dimensional model for $\eg$ and can be constructed to have either infinitely many conjugacy classes of finite subgroups, or to have centralisers of finite subgroups not of type $\FP_\infty.$
There are a number of classes of groups of type $\FP_\infty$ admitting cocompact models for $\eg$ including Gromov hyperbolic groups \cite{meintruppschick}, $Out(F_n)$ \cite{vogtmann} or elementary amenable groups of type $\FP_\infty$ \cite{kmn}.

\noindent Using our results on centralisers and conjugacy classes of finite subgroups we show:

\medskip\noindent{\bf Theorem \ref{Tmain}.} $T_r(\Sigma)$ is of type quasi-$\UFP_\infty$ if and only if $T_l(\Sigma)$ is of type $\FP_\infty$ for any $1\leq l\leq d$ such that $\text{gcd}(l,d)\mid r$.

\medskip\noindent We also consider the geometric analogue, to be of type quasi-$\underline\FFF_\infty,$ and the corresponding version of
\ref{Tmain}. We conjecture that similar results hold true for the groups $G_r(\Sigma)$.

\medskip\noindent The paper is structured as follows: In Section 2 we define the Cantor algebras and the corresponding generalisations of Thompson's groups $G$ and $T.$ We then use this Cantor algebra to build a model for $\eg$ in Section 3. In Section 4 we prove the results on centralisers and conjugacy classes of finite subgroups that will be used later.

In Section 5 we  collect all necessary background on Bredon cohomology with respect to an arbitrary family, and on Bredon cohomological finiteness conditions for modules. We prove an analogue to the Bieri-Eckmann criterion for property $\FP_n$ for modules. In Section 6 we specialise to the case of the family of finite subgroups and define what it means for a group to be quasi-$\UFP_\infty$ and quasi-$\UF_\infty.$ Finally, the main results are proven in Section 7.

\medskip\noindent{\bf Acknowledgements.} The authors wish to thank D.H. Kochloukova for very fruitful discussions,
without which, in fact, this work probably would not have happened. We also thank F. Matucci for a conversation, which led us to discover a gap in a previous version of this paper.

\section{Generalisations of Thompson-Higman groups}

\noindent As mentioned in the introduction, the generalised Thompson-Higman groups can be viewed as certain automorphisms groups of Cantor algebras. We shall begin by defining these algebras.
We  use the notation  of \cite{desiconbrita2}, Section 2. In particular, we
consider a finite set $\{1,\ldots,s\}$ whose elements are called colours. To each colour $i$ we associate an integer $n_i>1$ which is called its arity.
We say that $U$ is an $\Omega$-algebra, if, for each colour $i$, the following operations (we let all operations act on the right)
 are defined in $U$ (for detail, see \cite{Cohn} and \cite{desiconbrita2}):
 \begin{itemize}
\item[i)] One $n_i$-ary operation $\lambda_i$:
$$\lambda_i :U^{n_i}\to U.$$
We call these  operations ascending operations, or contractions.

\item[ii)] $n_i$ 1-ary operations $\alpha^1_i,\ldots,\alpha^{n_i}_i$:
$$\alpha^j_i:U\to U.$$
We call these operations 1-ary descending operations.

 \end{itemize}
 \noindent We denote $\Omega = \{ \lambda_i, \alpha_i^j \}_{i,j}$.
For each colour $i$ we also consider  the map $\alpha_i:U\to U^{n_i}$ given by
$$v\alpha_i:=(v\alpha^1_i,v\alpha^2_i,\ldots,v\alpha^{n_i}_i)$$
for any $v\in U$. These maps are called descending operations, or
expansions.
For any subset $Y$ of $U$, a simple expansion of colour $i$ of $Y$ is obtained by substituting some element $y\in Y$ by the $n_i$ elements of the tuple  $y\alpha_i$. A simple contraction of colour $i$ of $Y$ is the set obtained by substituting a certain collection of $n_i$ distinct elements of $Y$, say $\{a_1,\ldots,a_{n_i}\}$, by $(a_1,\ldots,a_{n_i})\lambda_i$. We also use the term operation to refer to the effect of a simple expansion, respectively contraction on a set .

For any set $X$ there is an  $\Omega$-algebra, free on $X$, which is called the  $\Omega$-word algebra on $X$ and is denoted by $W_\Omega(X)$. An admissible subset $A\subseteq W_\Omega(X)$ is a subset that can be obtained after finitely many expansions or contractions from the set $X$.

\noindent Descending operations can be visualised by tree diagrams, see the following example with $X=\{x\}$, $s=1$ and $n_1=2:$

\bigskip

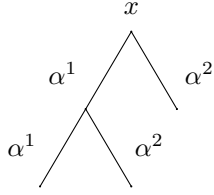
\begin{tikzpicture}[scale=0.6]

  \draw[black]
    (5,0) -- (6, 1.71) -- (7,0);
  \filldraw(6,1.71) circle (0.3 pt);

  \draw[black](4,-1.71)--(5,0)--(6, -1.71);

  \draw (5.5,0.8) node[left=8pt]{$\alpha^1$};
   \draw (6.5,0.8) node[right=8pt]{$\alpha^2$};
    \draw (5.4,-0.8) node[right=8pt]{$\alpha^2$};
     \draw (4.6,-0.8) node[left=8pt]{$\alpha^1$};
      \filldraw(5,0) circle (0.3pt);
        \filldraw(7,0) circle (0.3pt);
         \filldraw(6,-1.71) circle (0.3pt);
          \filldraw(4,-1.71) circle (0.3pt);
\draw(6, 1.71) node[above=3pt]{$x$};


   \end{tikzpicture}

The set $A=\{x\alpha^1\alpha^1,x\alpha^1\alpha^2,x\alpha^2\}$ is an admissible subset. In pictures we often omit the maps and label the nodes by positive integers as follows:

\bigskip

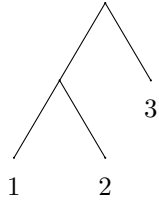
\begin{tikzpicture}[scale=0.6]

  \draw[black]
    (5,0) -- (6, 1.71) -- (7,0);
  \filldraw(6,1.71) circle (0.3 pt);

  \draw[black](4,-1.71)--(5,0)--(6, -1.71);

      \filldraw(5,0) circle (0.3pt);
        \filldraw(7,0) circle (0.3pt);
         \filldraw(6,-1.71) circle (0.3pt);
          \filldraw(4,-1.71) circle (0.3pt);

\draw(4,-1.71) node[below=4pt]{1};
\draw(6,-1.71) node[below=4pt]{2};
\draw(7,0) node[below=4pt]{3};

   \end{tikzpicture}

\noindent From now on  we  fix the set $X$ and assume it is finite. We  consider the variety of $\Omega$-algebras satisfying a certain set of identities as follows:

\begin{definition}\label{sigmadef} Let $\Sigma$ be the following set of
  laws in the alphabet $X$.

\begin{itemize}

\item[i)] For any $u\in W_\Omega(X)$, any colour $i$, and any $n_i$-tuple $(u_1,\ldots,u_{n_i})\in W_\Omega(X)^{n_i},$
 $$u\alpha_i\lambda_i=u,$$
 $$(u_1,\ldots,u_{n_i})\lambda_i\alpha_i=(u_1,\ldots,u_{n_i}).$$
 The set of all these relations is denoted $\Sigma_1$

\item[ii)] A certain set
$$\Sigma_2=\bigcup_{1\leq i<i'\leq s}\Sigma_2^{i,i'}$$
such that each $\Sigma_2^{i,i'}$ is either empty or  consists  of all the laws of the following form: Consider first $i$ and fix a map $f:\{1,\ldots,n_{i}\}\to\{1,\ldots,s\}$. For each $1\leq j\leq n_{i}$, we see  $\alpha_{i}^j\alpha_{f(j)}$ as a set of length 2 sequences of descending operations and  
let $\Lambda_{i}=\cup_{j=1}^{n_{i}}\alpha_{i}^j\alpha_{f(j)}$.  Do the same for $i'$ (with a corresponding map $f'$) to get $\Lambda_{i'}$ and now fix a bijection $\phi:\Lambda_{i}\to\Lambda_{i'}$.
Then $\Sigma_2^{i,i'}$ is the set of laws
$$u\nu=u\phi(\nu)\quad \nu\in \Lambda_{i},u\in W_\Omega(X).$$
(Note that by an abuse of notation we omit the $u\in W_\Omega(X)$ when we specify $\Sigma_2$ in the examples below).

\end{itemize}
\end{definition}

When factoring out the fully invariant congruence $\mathfrak{q}$ generated by $\Sigma$, we obtain an $\Omega$-algebra   $W_\Omega(X)/\mathfrak{q}$ satisfying the identities in $\Sigma$. For detail of the construction the reader is referred to \cite[Section 2]{desiconbrita2}.

\begin{definition} Let $r=|X|$ and $\Sigma$ as in Definition \ref{sigmadef}. Then the algebra $W_\Omega(X)/\mathfrak{q}=U_r(\Sigma)$ is called a Cantor-Algebra.
\end{definition}

\noindent Moreover,
there is an epimorphism of $\Omega$-algebras
$$\begin{aligned}W_\Omega(X)&\twoheadrightarrow U_r(\Sigma)\\
A&\mapsto\bar A.\\ \end{aligned}$$
As in \cite{desiconbrita2} we say that $\Sigma$ is {\sl valid}  if for any admissible $Y\subseteq W_\Omega(X)$, we have $|Y|=|\bar Y|$. This condition implies that $U_r(\Sigma)$ is a free object on $X$ in the class of those $\Omega$-algebras which satisfy the identities $\Sigma$ above.

\noindent If the set $\Sigma$ used to define $U_r(\Sigma)$ is valid, we also say that $U_r(\Sigma)$ is valid.

\begin{example}\label{higmanalgebra}
Higman \cite{higman} defined an algebra $V_{n,r}$ with $|X|=r$, $s=1$ and arity $n$ as above with $\Sigma_2$ being empty. This algebra, which we call Higman algebra, is used in the original construction of the Higman-Thompson-groups $G_{n,r}.$ For detail see also \cite{brown}.  In particular, these algebras are valid \cite[Section 2]{higman}.
\end{example}

\begin{example}\label{genhigman} Higman's construction for arity $n=2$ can be generalised as follows \cite[Section 2]{desiconbrita2}: Let $s\geq 1$ and $n_i=2$ for all $1\leq i \leq s.$ Hence we consider the set of $s$ colours $\{1,\ldots,s\}$, all of which have arity 2, together with the relations:
$\Sigma:=\Sigma_1\cup\Sigma_2$
with $$\Sigma_2:=\{\alpha_i^l\alpha_j^t=\alpha_j^t\alpha_i^l\mid
1\leq i \not= j\leq s; l,t=1,2 \}.$$
 Then $\Sigma$ is valid (see \cite{desiconbrita2} Lemma 2.9).

\noindent Furthermore one can also consider  $s$ colours,  all of arity $n_i=n,$ for all $1\leq i\leq s.$ Let
$$\Sigma_2:=\{\alpha_i^l\alpha_j^t=\alpha_j^t\alpha_i^l\mid
1\leq i \not= j\leq s; 1\leq l,t\leq n \}.$$
Using the same arguments as in \cite[Section 2]{desiconbrita2} one can show that the  $\Sigma$ obtained in this way is also valid.

We call the resulting Cantor algebras $U_r(\Sigma)$ Brin algebras.

\medskip\noindent The following tree-diagram visualises the relations in $\Sigma_2$. Here $r=1, s=2$ and $n=2.$ We express an expansion of colour $1$ with dotted lines and an expansion of colour $2$ by solid lines.  The leaves with the same label are identified.

\medskip
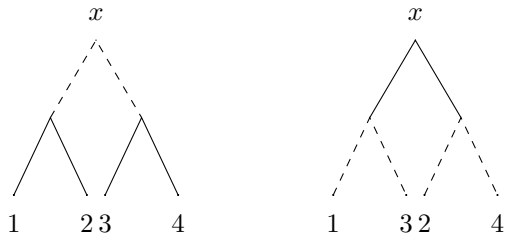
\begin{tikzpicture}[scale=0.6]

  \draw[black, dashed]
    (0,0) -- (1, 1.71) -- (2,0);
  \filldraw(1,1.71) circle (0.1pt) node[above=4pt]{$x$};

  \draw[black] (-0.8,-1.71) -- (0,0) --(0.8,-1.71);
  \draw[black] (1.2,-1.71) -- (2,0) --(2.8,-1.71);

      \filldraw (-0.8,-1.71) circle (0.3pt) node[below=4pt]{$1$};
      \filldraw (0.8,-1.71) circle (0.3pt) node[below=4pt]{$2$};
      \filldraw (1.2,-1.71) circle (0.3pt) node[below=4pt]{$3$};
      \filldraw (2.8,-1.71) circle (0.3pt) node[below=4pt]{$4$};

  \draw[black]
    (7,0) -- (8, 1.71) -- (9,0);
 \filldraw(8,1.71) circle (0.1pt) node[above=4pt]{$x$};

  \draw[black, dashed] (6.2,-1.71) -- (7,0) --(7.8,-1.71);
  \draw[black, dashed] (8.2,-1.71) -- (9,0) --(9.8,-1.71);

      \filldraw (6.2,-1.71) circle (0.3pt) node[below=4pt]{$1$};
      \filldraw (7.8,-1.71) circle (0.3pt) node[below=4pt]{$3$};
      \filldraw (8.2,-1.71) circle (0.3pt) node[below=4pt]{$2$};
      \filldraw (9.8,-1.71) circle (0.3pt) node[below=4pt]{$4$};

   \end{tikzpicture}

\end{example}

\begin{definition}
 Let $\Sigma$ be valid and consider $Y,Z\subseteq U_r(\Sigma)$. If $Z$ can be obtained from $Y$ by a finite number of simple expansions then we say that $Z$ is a {\sl descendant} of $Y$ and denote
$$Y\leq Z.$$
Conversely, $Y$ is called an {\sl ascendant} of $Z$ and can be obtained after a finite number of simple contractions. Note that this implies that if either of the sets $Y$ or $Z$ is admissible, then so is the other. In fact, the set of admissible subsets of $U_r(\Sigma)$ is a poset with respect to the partial order $\leq$. This poset is denoted by $\frak A_r(\Sigma)$.
 \end{definition}

It is easy to prove that any admissible subset is a basis of $U_r(\Sigma)$ (see \cite{desiconbrita2} Lemma 2.5).

\begin{remark} Let $\Sigma$ be valid and assume that we have $s$ colours of arities $\{n_1,\ldots,n_s\}$. Let $r$ be a positive integer. Observe that the cardinality of any admissible subset of $U_r(\Sigma)$ must be of the form  $m\equiv r$ mod $d$ for
$$d:=\text{gcd}\{n_i-1\mid i=1,\ldots,s\}.$$
Moreover, for any $m\equiv r$ mod $d,$ there is some admissible subset of cardinality $m$. And as admissible subsets are bases, we get $U_r(\Sigma)=U_m(\Sigma).$
\end{remark}

\begin{definition}\label{bd}  Let $B,C$ be admissible subsets of $U_r(\Sigma).$ We say that $T$ is the unique {\sl least upper bound} of $B$ and $C$ if $B \leq T$, $C \leq T$ and for all admissible sets $S$ such that $B\leq S$ and $C \leq S$ we have $T \leq S.$

\noindent We say, by abusing notation a little, that $U_r(\Sigma)$ is {\sl bounded}  if for all admissible subsets $B,C$ such that there is some admissible $A$ with $A\leq B,C$ there is a unique least upper bound of $B$ and $C$.
\end{definition}

\noindent One can also define greatest lower bounds, but this places a stronger restriction on the algebra, see \cite{desiconbrita2}. Moreover, note that a priori we require the existence of an upper bound only when our sets have a lower bound ($A$) but this turns out to be not too restrictive:

\begin{lemma}\label{upperbound} Let $U_r(\Sigma)$ be valid and bounded. Then any two admissible subsets have some (possibly not unique) common upper bound.
\end{lemma}
\begin{proof} Use the same proof as in  \cite[Proposition 3.4]{desiconbrita2}.
\end{proof}

\begin{example} The Brin algebras defined in Example \ref{genhigman} are valid and bounded. The existence of a unique least upper bound for $n=2$ is shown in \cite[Lemma 3.2]{desiconbrita2}. The general case is analogous.
\end{example}

\begin{example}\label{brownstein} Let $P\subseteq\Q_{>0}$ be a finitely generated multiplicative group. Then by a result of Brown, see  \cite[Proposition 1.1]{stein}, $P$ has a basis of the form $\{n_1,\ldots,n_s\}$ with all $n_i \geq 0$ ($i=1,...,s$).
Now consider $\Omega$-algebras on $s$ colours of arities $\{n_1,\ldots,n_s\}$ and let $\Sigma=\Sigma_1\cup\Sigma_2$ with $\Sigma_2$ the set of identities given by the following order preserving identification:
$$\{\alpha_i^1\alpha_j^1,\ldots,\alpha_i^1\alpha_j^{n_j},\alpha_i^2\alpha_j^1,\ldots,\alpha_i^2\alpha_j^{n_j},\ldots,\alpha_i^{n_i}\alpha_j^1,\ldots,\alpha_i^{n_i}\alpha_j^{n_j}\}=$$
$$\{\alpha_j^1\alpha_i^1,\ldots,\alpha_j^1\alpha_i^{n_i},\alpha_j^2\alpha_i^1,\ldots,\alpha_j^2\alpha_i^{n_i},\ldots,\alpha_j^{n_j}\alpha_i^1,\ldots,\alpha_j^{n_j}\alpha_i^{n_i}\},$$
where $i \neq j$ and $i,j\in \{1,...,s\}.$

\noindent The Cantor algebras $U_r(\Sigma)$ thus obtained will be called Brown-Stein algebras.
\end{example}

\noindent Note that, as $\{n_1,\ldots,n_s\}$ is a basis for $P$, the $n_i$ are all distinct. Hence, when visualising  the identities in $\Sigma_2$ for the Brown-Stein algebra, it suffices to only use one colour, as the arity of an expansion already determines the colour.
In the following example let $r=1, s=2, n_1=2$  and $n_2=3$.

\medskip
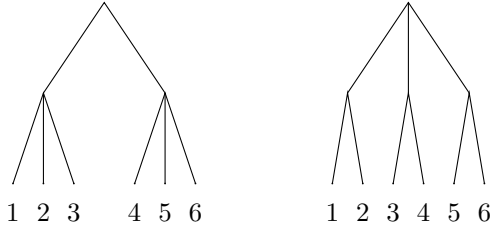
\begin{tikzpicture}[scale=0.4]

  \draw[black]
    (1,0) -- (3, 3) -- (5,0);

  \draw[black] (0,-3)--(1,0)--(2,-3);
  \draw[black] (1,-3)--(1,0);
  \draw[black] (4,-3)--(5,0)--(6,-3);
  \draw[black] (5,-3)--(5,0);

      \filldraw (0,-3) circle (0.3pt) node[below=4pt]{$1$};
      \filldraw (1,-3) circle (0.3pt) node[below=4pt]{$2$};
      \filldraw (2,-3) circle (0.3pt) node[below=4pt]{$3$};
      \filldraw (4,-3) circle (0.3pt) node[below=4pt]{$4$};
      \filldraw (5,-3) circle (0.3pt) node[below=4pt]{$5$};
      \filldraw (6,-3) circle (0.3pt) node[below=4pt]{$6$};

  \draw[black]  (11,0) -- (13, 3) -- (15,0);
  \draw[black] (13,0)--(13,3);

  \draw[black] (10.5,-3)--(11,0)--(11.5,-3);
 \draw[black] (12.5,-3)--(13,0)--(13.5,-3);
  \draw[black] (14.5,-3)--(15,0)--(15.5,-3);

      \filldraw (10.5,-3) circle (0.3pt) node[below=4pt]{$1$};
      \filldraw (11.5,-3) circle (0.3pt) node[below=4pt]{$2$};
      \filldraw (12.5,-3) circle (0.3pt) node[below=4pt]{$3$};
      \filldraw (13.5,-3) circle (0.3pt) node[below=4pt]{$4$};
      \filldraw (14.5,-3) circle (0.3pt) node[below=4pt]{$5$};
      \filldraw (15.5,-3) circle (0.3pt) node[below=4pt]{$6$};

\end{tikzpicture}

\begin{lemma}\label{validbrownstein} The Brown-Stein algebras are valid and bounded.
\end{lemma}
\begin{proof} This is Proposition 1.2 (due to K. Brown) in \cite{stein}. \end{proof}

\noindent In fact, in \cite{stein} Lemma \ref{upperbound} is proven directly, i.e that any two admissible subsets have some common upper bound.

\bigskip\noindent We can now define the generalised Thompson-Higman groups. Recall, that in a valid Cantor algebra $U_r(\Sigma)$, admissible subsets are bases.

\begin{definition}\label{groups} Let $U_r(\Sigma)$ be a valid Cantor algebra. We define  $G_r(\Sigma)$ to be the group of those $\Omega$-algebra automorphisms of $U_r(\Sigma)$, which are induced by a  map $V\to W$ , where $V$ and $W$ are admissible subsets of the same cardinality.
\end{definition}

\begin{example} If $U_r(\Sigma)$ is a Higman algebra as in Example \ref{higmanalgebra}, we retrieve the original Higman-Thompson-groups $G_{n,r}$.
Let $U_r(\Sigma)$ be a Brin algebra on $s$ colours of arity $2$ as in Example \ref{genhigman}. Then the groups constructed are Brin's \cite{Brin} generalisations $sV$ of Thompson's group $V=G_{2,1}.$ The description of $sV$ as automorphism groups of a Cantor algebra can be found in \cite{desiconbrita2}. Finally, the groups $G_r(\Sigma)$, when $U_r(\Sigma)$ is a Brown-Stein algebra as in Example \ref{brownstein}, were considered in \cite{stein}.
\end{example}

\begin{remark} It is conceivable that in fact $G_r(\Sigma)$ equals the full group of  $\Omega$-algebra automorphisms of $U_r(\Sigma)$. This would follow if one could prove that any finite basis of $U_r(\Sigma)$ is an admissible subset (this is the case for the Higman algebra, see \cite{higman} Corollary 1). 
\end{remark}

\bigskip\noindent We go back to the case of an arbitrary valid Cantor algebra $U_r(\Sigma)$ and assume that the set $X$ is ordered. It can be seen that this order is  inherited by certain subsets of $W_\Omega(X)$ including all admissible subsets, see for example \cite{brown} or \cite{higman}. If the relations in $\Sigma_2$ preserve that ordering, in the sense that the bijection $\phi$ in Definition \ref{sigmadef} do, then we also have an inherited order on the admissible subsets of $U_r(\Sigma)$. We shall call this the induced ordering.

\begin{definition}
Suppose we have a Cantor algebra $U_r(\Sigma)$ where $\Sigma$ preserves the induced ordering. We may define  subgroups $F_r(\Sigma)$ and $T_r(\Sigma)$  of $G_r(\Sigma)$ as follows. We let $F_r(\Sigma)$ be the group of order preserving automorphisms between ordered admissible subsets of the same cardinality and $T_r(\Sigma)$ the group of cyclic order preserving automorphisms between ordered admissible subsets of the same cardinality.
\end{definition}

\begin{example}\label{ftexamples}  For $U_r(\Sigma)$ a Higman algebra of Example \ref{higmanalgebra} the definition above yields the groups $F_{n,r}$ and $T_{n,r}$ as in \cite{brown}. Recall that Thompson's groups are $F=F_{2,1}$ and $T=T_{2,1}$.

\noindent
Let $U_r(\Sigma)$ be a Brown-Stein algebra as in Example \ref{brownstein} In this case, $\Sigma$ is order preserving, so we may define the groups $F_r(\Sigma)$ and $T_r(\Sigma)$, which are considered in \cite{stein}.

\noindent Since $\Sigma_2$ in the definition of the Brin algebra of Example \ref{genhigman} is not order-preserving, it there is no obvious way to define the groups $F_r(\Sigma)$ or $T_r(\Sigma)$ for this algebra.
\end{example}

\begin{remark} Note, that if definable, the groups $F_r(\Sigma)$ are torsion-free. In both cases mentioned in Example \ref{ftexamples}, the resulting groups $F_r(\Sigma)$ are known to be of type $\FP_\infty$ and finitely presented \cite{brown, stein}.

\noindent Since, for torsion-free groups ordinary and Bredon cohomological finiteness conditions are identical, we will not consider these groups further.
\end{remark}

\section{A model for $\eg$ for generalised Thompson groups}

\noindent From now on we fix a valid  $\Sigma$ and a finite positive integer $r$. Also assume that the Cantor algebra $U_r(\Sigma)$ is bounded. In this section we give a quite elementary proof of the following result.

\begin{theorem}\label{model} The geometric realisation of the poset of admissible subsets is a model for $\underbar{E}G_r(\Sigma)$.
\end{theorem}


We fix an admissible subset $X\subseteq U_r(\Sigma)$ of cardinality $r$.

\begin{lemma}\label{fin1} For any finite $Q\leq G_r(\Sigma)$ there exists some admissible subset $Z$ such that $ZQ=Z$. Moreover we may assume $X\leq Z$.
\end{lemma}
\begin{proof} For every $q\in Q$ choose a common upper bound $T_q$ of $X$ and $Xq$. Then put $Z_q:=T_qq^{-1}$ and let $Y$ be an upper bound of
$$\{Z_q\mid q\in Q\}.$$
Note that $X\leq Z_1=T_1$ and for any $q\in Q$,
$$X\leq T_q=Z_qq\leq Yq.$$
Therefore we may choose $Z$ the  least upper bound of $\{Yq\mid q\in Q\}$.  By definition of unique least upper bound we get $ZQ=Z$.
\end{proof}

\begin{proposition}\label{fin2} Any two elements in $\frak A_r(\Sigma)^Q$ have an upper bound in $\frak A_r(\Sigma)^Q.$
\end{proposition}
\begin{proof} Let $Y,Z\in\frak A_r(\Sigma)^Q$. We begin by showing that there are admissible sets $Y_1,Z_1\in\frak A_r(\Sigma)^Q$ such that $Y_1$ is an upper bound of $X$ and $Y$ and $Z_1$ is an upper bound of $X$ and $Z$. It suffices to prove the existence of $Y_1$.
Take  an upper bound $Y_2\in \frak A_r(\Sigma)$ of $X$ and $Y$ and consider
$$\{Y_2q^{-1}\mid q\in Q\}.$$
Let $Y_3\in \frak A_r(\Sigma)$ be an upper bound  of this set. Then, for any $q\in Q,$
$$Y_2\leq Y_3q.$$
 Therefore $X\leq Y_3q$. This implies that we may choose $Y_1$ to be the least upper bound of
$$\{Y_3q\mid q\in Q\}.$$
Clearly, $Y,X\leq Y_1$. Again, the definition of least upper bound  implies that $Y_1\in\frak A_r(\Sigma)^Q$.

Now, let $T$ be the least upper bound of $Y_1$ and $Z_1$. Then for any $q\in Q$
 $$Y_1=Y_1q\leq Tq,$$
 $$Z_1=Z_1q\leq Tq$$
 so we get $T\in\frak A_r(\Sigma)^Q.$
 \end{proof}

\begin{proof}(of Theorem \ref{model})
 Lemmas \ref{fin1} and \ref{fin2} imply that for any finite subgroup $Q\leq G_r(\Sigma)$ the poset $\frak A_r(\Sigma)^Q$ is non-empty and directed, thus $|\frak A_r(\Sigma)|^Q=|\frak A_r(\Sigma)^Q|\simeq\ast$. Moreover for any $V\in\frak A_r(\Sigma)$, $$\text{Stab}_{G_r(\Sigma)}(V)=\{g\in G_r(\Sigma)\mid Vg=V\}$$
 is contained in the group of permutations of the finite set $V$, thus it is finite. This implies that for any
 $H\leq G_r(\Sigma)$, $\frak A_r(\Sigma)^H=\emptyset$ unless $H$ is finite.
\end{proof}

\noindent
This model is not of finite type, but there is a filtration of $|\frak A_r(\Sigma)^Q|$ by finite type subcomplexes, exactly as in the construction in \cite[Theorem 4.17]{brown}:

\begin{proposition}\label{filtration} For any finite $Q\leq G_r(\Sigma)$ there is a filtration of $|\frak A_r(\Sigma)^Q|$
$$\ldots\subset|\frak A_r(\Sigma)^Q|_{h-1}\subset|\frak A_r(\Sigma)^Q|_h\subset|\frak A_r(\Sigma)^Q|_{h+1}\subset\ldots$$ such that each $|\frak A_r(\Sigma)^Q|_h/C_{G_r(\Sigma)}(Q)$ is finite.\end{proposition}
\begin{proof} Let
$$|\frak A_r(\Sigma)^Q|_h:=\{Y\in\frak A_r(\Sigma)^Q\mid|Y|\leq h\}.$$
Consider $Y,Z\in\frak A_r(\Sigma)^Q$ with $|Y|=|Z|$ and isomorphic as $Q$-sets. This means that there is a $Q$-bijection
$$\sigma:Y\to Z.$$
Let $g\in G_r(\Sigma)$ be the element given by $yg=y\sigma$ for each $y\in Y$. Then for any $q\in Q$, $(yq)g=(yq)\sigma=y\sigma q=ygq$. This means that the commutator $[g,q]$
acts as the identity on the admissible set $Y$ and therefore $[g,q]=1$. Hence $g\in C_{G_r(\Sigma)}(Q)$. As for any $m\leq h$ there are finitely many possible $Q$-sets of cardinality $m$, the result follows.
\end{proof}

\begin{remark} Provided that $\Sigma$ is order-preserving, Theorem \ref{model} and Proposition \ref{filtration} can be restated replacing $G_r(\Sigma)$ with $T_r(\Sigma).$

\end{remark}

\begin{remark}  The filtration of Proposition \ref{filtration} is used by Brown \cite[Theorem 4.17]{brown} with $Q=1,$ to show that the Higman-Thompson groups $G_{n,r}, T_{n,r}$ and $F_{n,r}$  are of type $\FP_\infty$. The approach used by Brown is as follows: Fix an admissible subset $Y$. Show that if $|\frak A_{n,r}|_{<Y}$ denotes the set of admissible subsets which are contractions of $Y$, then the connectivity of $|\frak A_{n,r}|_{<Y}$  grows with the cardinality of $Y$. Then, show that this implies that the connectivity of the pair $(|\frak A_{r,n}|_{h+1},|\frak A_{r,n}|_{h+1})$ tends to $\infty$, which in turn yields that $G_{n,r}, T_{n,r}$ and $F_{n,r}$ are all of type $\FP_\infty$. Key to this approach is understanding the complex $|\frak A_{n,r}|_{<Y}$. In the case of the Higman algebra, Brown shows (\cite[Lemmas 4.18; 4.19]{brown}), that any two simple contractions $Y_1,Y_2$ of $Y$ have a common lower bound if and only if the contracted vertices are disjoint, which allows him to show that $|\frak A_{n,r}|_{<Y}$ is homotopy equivalent to a much simpler complex. However, this is no longer true if we work with a more general Cantor algebra $U_r(\Sigma)$: Consider for example a Brown-Stein algebra as in Example \ref{brownstein} with  arities 2 and 3. Let $Y$ be any admissible set with 6 elements labeled 1, 2, 3, 4, 5 and 6. Let $Y_1$ be the simple contraction of arity 2 of the elements 3 and 4 and $Y_2$  the simple contraction of arity 3 of the elements 1, 2 and 3. Then the sets of contracted vertices are not disjoint, however there is a common lower bound $Z\leq Y_1,Y_2$ as the picture before Lemma \ref{validbrownstein} shows.  Stein  used a different method to the one described here to prove that the groups of \cite {stein} are of type $\FP_\infty.$

Similar problems were encountered when Kochloukova and the authors considered Brin's groups  \cite{desiconbrita2}. In general, the same difficulty applies to the groups $G_r(\Sigma)$, as well as to  $T_r(\Sigma), F_r(\Sigma)$ where definable. It is conceivable, however, that  Brown's approach can be applied more generally using  an analogue of Brown's connectivity result, see for example \cite{desiconbrita2} where it is used to show that  Brin's groups for $r=1$, $s=2,3$ are of type $\FP_\infty.$
\end{remark}

\section{Centralisers and conjugacy classes of finite subgroups for $G_r(\Sigma)$ and $T_r(\Sigma).$}

\noindent Let $Q\leq G_r(\Sigma)$ be a finite subgroup. In this section we give a more detailed analysis of the poset $\frak A_r(\Sigma)^Q$ to describe $C_{G_r(\Sigma)}(Q)$ and the number of conjugacy classes of subgroups isomorphic to $Q$.  In case  $T_r(\Sigma)$ is defined, we also derive the corresponding results.  This will be used later when we  prove our main result on the  cohomological finiteness properties of these groups.

 Let $\{w_1,\ldots,w_t\}$ be the set of lengths of all the possible transitive permutation representations of $Q$.
Any $Y\in\frak A_r(\Sigma)^Q$ is a finite $Q$-set so it is determined by its decomposition into transitive $Q$-sets. If we take one of those sets and apply the operation $\alpha_i$ for a fixed colour $i$ to each of its elements, we obtain a new admissible subset which is also fixed by $Q$. We say that this is a simple $Q$-expansion of $Y$. More explicitly, the admissible set obtained from $Y$ is:
$$Y\setminus\{yq\mid q\in Q\}\cup\{yq\alpha_i^j\mid q\in Q,1\leq j\leq n_i\}$$
for a certain $y\in Y$. We also use the term $Q$-expansion to refer to a chain of simple $Q$-expansions.

Conversely, if we choose $n_i$ different orbits of the same type (i.e., corresponding to the same permutation representation)  in $Y$, then we may contract them to a single orbit (of the same type). We call this a simple $Q$-contraction. Simple $Q$-contractions are more complicated to handle than simple $Q$-expansions: we may contract an element of the first of the orbits with any of the elements on the others. Hence, even if the orbits to be contracted are determined, there are many possibilities to perform the explicit contraction.
Note that the admissible subsets obtained this way will lie in $\frak A_r(\Sigma)^Q$.

Large parts of the next three results can be found in \cite[Section 6]{higman}. We shall, for the reader's convenience,  recall  the arguments.

\begin{lemma}\label{fin3} Let $Y,Z\in\frak A_r(\Sigma)^Q$ with $Y<Z$ and assume there is no admissible subset $C\in\frak A_r(\Sigma)^Q$ with $Y\leq C\leq Z$. Then $Z$ is a simple $Q$-expansion of $Y$. Hence $Y$ is a simple $Q$-contraction of $Z$.
\end{lemma}
\begin{proof} We may choose  a chain of simple expansions
$$Y<Y_1<\ldots<Y_r<Z.$$
Let $w\in Y$ be the vertex expanded in the first simple expansion $Y<Y_1$ and $W\subseteq Y$ be the $Q$-orbit with $w\in W$. Assume also that this first expansion corresponds to the colour $i$. Then as $Z$ contains certain descendants of $\{w\alpha_i\}$ and it is $Q$-invariant it must also contain the analogous descendants of $\{u\alpha_i\mid u\in W\}$. Therefore if $C$ denotes the simple $Q$-expansion consisting of expanding $W$ by $\alpha_i$, then $Y<C\leq Z$. As $C\in\frak A_r(\Sigma)^Q$, we deduce by the hypothesis that $C=Z$.
\end{proof}

\begin{proposition}\label{fixed1} For any finite subgroup $Q\leq G_r(\Sigma)$, there is a uniquely determined set of integers $\pi(Q):=\{r_1,\ldots,r_t\}$ with $0\leq r_j\leq d$ and
$$\sum_{j=1}^tr_jw_j\equiv r\text{ mod }d$$
such that there is an admissible subset $Y\in\frak A_r(\Sigma)^Q$ with $|Y|=\sum_{j=1}^tr_jw_j$.

\noindent Moreover, any other element in $\frak A_r(\Sigma)^Q$ can be obtained from $Y$ by a finite sequence of simple $Q$-expansions or $Q$-contractions.
\end{proposition}
\begin{proof} First, note that by \ref{fin1}, $\frak A_r(\Sigma)^Q\neq\emptyset$. Now choose some $Z\in\frak A_r(\Sigma)^Q$ and decompose it as a disjoint union of transitive $Q$-sets. Let $k_j$ be the number of transitive sets in this decomposition which are of type $j$, i.e which correspond to the same permutation representation.  Observe that whenever we apply simple $Q$-contractions or $Q$-expansions to $Z$, if the set thus obtained has $m_j$ transitive $Q$-sets of type $j$, then $m_j\equiv k_j$ mod $d$. Note also that
$$|Z|=\sum_{j=1}^tk_jw_j\equiv r\text{ mod }d.$$
Let
$$r_j=\Bigg\{
\begin{aligned}
&0,\text{ if }k_j=0\\
&d,\text{ if }0\neq k_j\equiv 0\text{ mod }d\\
&l\text{ with }0<l<d\text{ and }l\equiv k_j\text{ mod }d,\text{ otherwise.}\\
\end{aligned}
$$
 By successively performing simple $Q$-contractions or $Q$-expansions of $Z$ we may get an admissible set $Y$ such that the number of transitive $Q$-sets of type $j$ in $Y$ is exactly $r_j$. Observe that the $r_j$ are uniquely determined, whereas $Y$ is not.
Finally, \ref{fin2} implies that for any other $C\in\frak A_r(\Sigma)^Q$, there is an upper bound, say $D$, of $Y$ and $C$ with $D\in\frak A_r(\Sigma)^Q$ which means that
$$Y\leq D\geq C.$$
 By Lemma \ref{fin3} we may choose chains
 $$Y=D_0<D_1<\ldots<D_{l_1}=D=C_0> C_1>\ldots >C_{l_2}=C$$
 such that each step consists of a simple $Q$-expansion/$Q$-contraction and we are done.
\end{proof}

\begin{theorem}\label{conjclasses} Let $Q_1,Q_2\leq G_r(\Sigma)$ be finite subgroups with $Q_1\cong Q_2$. Then $Q_1$ and $Q_2$ are conjugate in $G_r(\Sigma)$ if and only if $\pi(Q_1)=\pi(Q_2)$.

\noindent In particular, there are only finitely many conjugacy classes of subgroups isomorphic to $Q_1.$
\end{theorem}

\begin{proof}  Fix an isomorphism $\phi:Q_1\to Q_2$. Assume first that $\pi(Q_1)=\pi(Q_2)$. Then there are admissible subsets $V_1,V_2$ with $V_i\in\frak A_r(\Sigma)^{Q_i}$ having the same number of elements and moreover the same structure as $Q_i$-sets, which means that there is a bijection between them which we denote $g$ such that for any $q\in Q_1$ and $v\in V_1$, $(vq)g=vg q^\phi$. This yields an element $g\in G_r(\Sigma)$ with $g^{-1}qg=q^\phi$.

Conversely, assume $Q_2=g^{-1}Q_1g$ with $g\in G_r(\Sigma)$. Then for any $V_1\in\frak A_r(\Sigma)^{Q_1}$, $V_1g\in\frak A_r(\Sigma)^{Q_2}$. Moreover, $g$ induces an isomorphism as $Q_i$-sets so the orbit structure of the minimal elements of $\frak A_r(\Sigma)^{Q_1}$ and $\frak A_r(\Sigma)^{Q_1}$ has to be the same.

\end{proof}

\begin{theorem}\label{centralizer} Let $Q\leq G_r(\Sigma)$ be a finite subgroup and $\pi(Q)=\{r_1,\ldots,r_t\}$ as in Proposition \ref{fixed1}. Then
$C_{G_r(\Sigma)}(Q)\cong H_{r_1}\times\ldots\times H_{r_t}$ where each of the $H_{r_i}$ fits into the following split group extension 
$$K_{i} \mono H_{r_i} \twoheadrightarrow G_{r_i}(\Sigma)$$
with $K_i$ locally finite.
\end{theorem}

\begin{proof} Choose an admissible $Y\in\frak A_{r}(\Sigma)^Q$ as in Proposition \ref{fixed1}. We begin by proving the result in the special case when there are exactly $l$ $Q$-orbits all of the same type in the $Q$-set $Y.$ In other words, we  assume that  in  Proposition \ref{fixed1} for some $k$, $l:=r_k$ and all the others $r_j=0$. 
Let $w:=w_k$ be the length of those $Q$-orbits, and for each $i=1,...,l$, choose an orbit  representative $y_i.$ 
We call the subset $\{y_1,...,y_l\}\subseteq Y$ thus obtained the set of marked elements. 
We consider any $Q$-expansion of $Y$ as marked, by marking precisely the descendants of marked elements in $Y$. And we say that a $Q$-contraction is marked if marked elements are contracted only with marked elements and result in the marked elements of the new subset. Note that this implies that elements of the form $yq$ with $q\in Q$ and $y$ marked can only be contracted with elements $y'q$ for the same $q\in Q$ and $y'$ marked. 
We now define: 
$$\begin{array}{ll}\mathcal {M}= & \{ M \,|\, M \in \frak A_r(\Sigma)^Q \mbox{admissible and obtained from $Y$ by marked}  \\ &\mbox{$Q$-expansions and marked $Q$-contractions}\}.\end{array}$$
$\mathcal M$ is the set of marked admissible subsets of $\frak A_r(\Sigma)^Q$ and can also be seen as the diagonal subposet:
$$\mathcal{M} \subseteq {\underbrace{\frak A_l(\Sigma)\times \ldots\times\frak A_l(\Sigma)}_w}.$$


Now fix an admissible subset $X=\{x_1,\ldots,x_l\}\in\frak A_l(\Sigma)$ and fix a bijection $\iota_{X,Y}:x_i\mapsto y_i$ from $X$ to the marked elements of $Y$. From $X\mapsto Y$ we get a poset map
$$\iota^{\frak A}:\frak A_l(\Sigma)\to\frak A_r(\Sigma)^Q,$$
which commutes with expansions and contractions in $\frak A_l(\Sigma)$, and with $Q$-expansions and marked $Q$-contractions in $\frak A_r(\Sigma)^Q$. This is well defined since $X$ is a basis of the algebra used to construct $\frak A_l(\Sigma)$. The fact that we only contract marked elements allows us to avoid ambiguities. Observe that  $\text{Im}(\iota^{\frak A})=\mathcal{M}$.
 Moreover, whenever $\iota^{\frak A}(X_1)=Y_1$, there is a well defined bijection $\iota_{X_1,Y_1}$ between $X_1$ and the set of marked elements in $Y_1.$ For convenience we let $\iota_{X,Y}$ act  on the right. We use this to define a group homomorphism 
 $$\iota:G_l(\Sigma)\to C_{G_r(\Sigma)}(Q),$$
as follows: Let $g$ be given by a map $g: X_1\to X_2$ for $X_1,X_2\in\frak A_r(\Sigma)$ and put $Y_1=\iota^{\frak A}(X_1)$, $Y_2=\iota^{\frak A}(X_2)$. Then  $\iota(g):Y_1\to Y_2$ is the unique map which commutes with the $Q$-action and such that $g\iota_{X_2,Y_2}=\iota_{X_1,Y_1}\iota(g)$ (recall that the marked elements are representatives of the $Q$-orbits).
Obviously $\iota(g)\in C_{G_r(\Sigma)}(Q)$.

\noindent Next, we define a second poset map
$$\tau^{\frak A}:\frak A_r(\Sigma)^Q\to\frak A_l(\Sigma)$$
such that $\tau^{\frak A}\iota^{\frak A}=id_{\frak A_l(\Sigma)}$. To do this, put $\tau^{\frak A}(Y)=X$, identify all the elements in the $Q$-orbit of each $y_i$ with $x_i$  and extend using the corresponding operations on both sides. Proposition \ref{fixed1} and the fact that $Y$ is admissible, imply that $\frak A_r(\Sigma)^Q$ is also free on $Y$, hence $\tau^{\frak A}$ is well defined. In an analogous way as before, there is also an explicit bijection between the $Q$-orbits in any $Y_1$ and the elements of $\tau^{\frak A}(Y_1)$ which can be used to define a group homomorphism
$$\tau: C_{G_r(\Sigma)}(Q)\to G_{l}(\Sigma).$$
Observe that whenever $g\in C_{G_r(\Sigma)}(Q)$ and $Y_1\in\frak A_r(\Sigma)^Q$, then $Y_1g\in\frak A_r(\Sigma)^Q$.

In particular, $\tau\iota = id_{G_l(\Sigma)}, $ giving us the desired split group extension.  We now proceed to describe $K:=\text{Ker}\tau$. To begin we observe that
$K$ consists precisely of those $h\in C_{G_r(\Sigma)}(Q)$ such that for any $A\in \frak A_r(\Sigma)^Q$, $\tau^{\frak A}(Ah)=\tau^{\frak A}(A)$ and $h$ fixes the $Q$-orbits of $A$ setwise. 

We claim that for any $h\in K$ there is some $Q$-expansion of $Y$, $Z\in\frak A_r(\Sigma)^Q$ with $Zh=Z$. To see this, using Proposition \ref{fin2}, take $Z\in\frak A_r(\Sigma)^Q$ to be an upper bound of $Y$, $Yh^{-1}$. Then Lemma \ref{fin3} implies that $Z$ and $Zh$ are both $Q$-expansions of $Y$ and therefore they are marked. Thus $Z,Zh\in\mathcal{M}=\text{Im}\iota^{\frak A}$. As $h$ lies in $K$, we have $\tau^{\frak A}(Z)=\tau^{\frak A}(Zh)$. So the fact that $\tau^{\frak A}$ is injective when restricted to $\text{Im}\iota^{\frak A}$ implies the claim.
In particular, $K$ is the union of its subgroups of the form
$$K_Z:=\{h\in C_{G_r(\Sigma)}(Q)\mid Zh=Z\text{, $h$ fixes the $Q$-orbits setwise}\}$$
where $Z$ is a $Q$-expansion of $Y$.  As each $K_Z$ is finite, using Proposition \ref{fin2} we see that $K$ is locally finite, thus proving the special case.

To finish our proof, we now prove the general case when $Y$ has $Q$ orbits of different types. Let $\pi(Q)=\{r_1,\ldots,r_t\}$ and $w_1,\ldots,w_t$ be as in Proposition \ref{fixed1}. Let $Y=\bigcup_{i=1}^tY_i$ with $Y_i$ the union of the $r_i$ $Q$-orbits of type $i$ in $Y$. Then $Q$ acts on each $Y_i.$ Note that a single action might not be faithful, but the intersection of the kernels must be trivial.
Also note, that in $\frak A_r(\Sigma)^Q$, $Q$-contractions can not mix elements belonging to orbits of different type. This implies that we have a direct product of posets
$$\frak A_r(\Sigma)^Q\cong\frak A_{w_1r_1}(\Sigma)^{Q}\times\ldots\times\frak A_{w_tr_t}(\Sigma)^Q,$$
where we let the group $Q$ act on each poset $\frak A_{w_ir_i}(\Sigma)$ using its action on $Y_i$ and extending via extensions and contractions. This action yields also a group homomorphism $\phi_i:Q\to G_{w_ir_i}(\Sigma).$ The direct product of posets above implies that $C_{G_r(\Sigma)}(Q)$ decomposes as the direct product of the centralisers of $\phi_i(Q)\leq G_{w_ir_i}(\Sigma)$. For each of these we can apply the case of  a single type of orbit, and we are done.
\end{proof}

\begin{remark} In an analogous way, one can prove that there is also a group epimorphism
$$N_{G_r(\Sigma)}(Q)\twoheadrightarrow  G_{r_1}(\Sigma)\times\ldots\times G_{r_t}(\Sigma)$$
with locally finite kernel. 
\end{remark}

\begin{remark} In \cite{matuccietc}, there is a description of centralisers of elements $g$ in the Higman groups $G_{n,1}$ associated to $U_r(\Sigma)=V_{n,1}.$ Whenever $Q=<g>$ has finite order in those groups, this coincides with ours.
\end{remark}

\begin{remark}\label{kernel}
With little more effort we can give a description of  the kernel $K$ appearing in the single type of orbit case in the proof of Theorem \ref{centralizer}: Let $S_w$ be the symmetric group of degree $w$ and choose a permutation representation $\phi:Q\to S_w$ associated to the $Q$ action on the orbits of $Y$. Denote $L:=C_{S_{w}}(\phi(Q))$. Then if $Z$ is a $Q$-expansion of $Y$ with $|Z|=mw$, we have the following isomorphism:
$$\{h\in C_{G_r(\Sigma)}(Q)\mid Zh=Z\text{, $h$ fixes the $Q$-orbits setwise}\}\cong L^m:=\underbrace{L\times\ldots\times L}_m.$$

Consider now a simple $Q$-expansion $Z\leq Z_1$ with $|Z_1|=m_1w,$ consisting of applying a descending operation of arity $n_j$ to a vertex $z_0.$ We get a group homomorphism $L^m \to L^{m_1}$ given by mapping the copy of $L$ corresponding to $z_0$ to the product of copies of $L$ corresponding to the descendants of $z_0$ via the diagonal map $L\to L^{n_j}$, and leaving the remaining factors intact. This gives, in an obvious way, a direct system of groups and hence  $K$ is the directed limit of the system thus obtained. 

\end{remark}

\medskip\noindent We shall now consider the groups $T_r(\Sigma)$ whenever they are defined, i.e. whenever $\Sigma$ preserves the induced order in the admissible subsets of our Cantor-Algebra $U_r(\Sigma)$. In this case, centralisers of finite subgroups have an easier structure.

  A first observation is that any finite subgroup $Q \leq T_r(\Sigma)$ is cyclic.  Moreover following the argument of Proposition \ref{fixed1} we see that by writing the transitive permutation representations of $Q$ to have the faithful representation first, i.e. $w_1=|Q|$, we obtain  $\pi(Q)=\{r_1,0,\ldots,0)$ and $r_1w_1\equiv r\text{ mod }d.$ To see this, take for example $Y\in\frak A_r(\Sigma)^Q$ the admissible subset obtained following the argument of Proposition \ref{fixed1} and assume that certain $g\in Q$ fixes some $y_0\in Y$. The condition that $g$ preserves cyclically the order, implies that $g$ fixes $Y$ pointwise, thus $g=1$. As a consequence, $|Y|=r_1w_1$.

\begin{theorem}\label{Textension}
 Let $U_r(\Sigma)$ be a Cantor-Algebra with order preserving $\Sigma$ and
 $Q \leq T_r(\Sigma)$ a finite subgroup. Then there is only one conjugacy class of finite subgroups of order $|Q|$ and for a certain $0<l\leq d,$ depending on $Q$ there is a central extension
$$Q \mono C_{T_r(\Sigma)}(Q) \epi T_l(\Sigma).$$

\end{theorem}

\begin{proof} For the first assertion observe that any two cyclic groups of the same order are isomorphic and they only have one faithful permutation representation. 
Hence it suffices to choose a cyclic order preserving $h$ between the corresponding admissible subsets $Y$. Note that they have the same cardinality. 

For the second assertion, embed $T_r(\Sigma)$ in $G_r(\Sigma)$ and let $\iota,\tau$ be the group homomorphisms of Theorem \ref{centralizer}, we use the same notation as there. The result will follow once we check that $\tau(C_{T_r(\Sigma)}(Q))=T_l(\Sigma)$ and that $Q=\text{Ker}\tau\cap T_r(\Sigma)$. 

Note first that we may choose the map $\iota_{X,Y}$ to be order preserving. 
The fact that the action of $Q$ cyclically preserves  the order on $Y,$ implies that we may assume that if for any basis $Y_1$ the marked elements are $\{y_1,\ldots,y_m\}$, then the elements of $Y_1$ are ordered as 
$$y_1<\ldots<y_m<y_1q<\ldots<y_mq<\ldots<y_1q^i<\ldots y_mq^i<\ldots,$$
for certain (fixed) $q$ generating $Q$ (note that the marked elements of elements of $Y$ can be chosen so that $Y$ is ordered this way).
If $g\in C_{T_r(\Sigma)}(Q)$ represents a map between two such sets cyclically preserving that ordering, then it is obvious that the corresponding map $\tau(g)$ also does. In fact, if we denote  $y'_1<\ldots<y'_m<y'_1q<\ldots<y'_mq<\ldots<y'_1q^i<\ldots y'_mq^i<\ldots$
the elements of $Y_1g$ and choose the index $j$ such that $y_jg=y'_mq^a$. Then if $j<m$, $y_{j+1}g=y'_1q^{a+1}$ and if $j=m$,  $y_1g=y'_1q^a$. Here $0\leq a\leq |Q|-1$. This implies $\tau(C_{T_r(\Sigma)}(Q))\subseteq T_l(\Sigma)$. 

Conversely, take $g\in T_l(\Sigma)$. Then $g$ is determined by its action on a pair of ordered admissible subsets $X_1:x_1<\ldots<x_m$, $X_2=X_1g:x'_1<\ldots<x'_m$. Put $Y_1:=\iota^{\frak A}(X_1)$, $Y_2:=\iota^{\frak A}(X_2)$ and denote their elements as before. Now, $\iota(g)$ as defined in Theorem \ref{centralizer} does not cyclically preserve  the order between $Y_1$ and $Y_2$. Let $j$ be the subindex such that $x_1g=x'_j$. We construct  $k\in K=\text{Ker}\tau$ as follows: 
$$y'_iq^ak=\Bigg\{\begin{aligned}
&y'_iq^{a+1}\text{ for }1\leq i<j\\
&y'_iq^a\text{ for }j\leq i\leq m,\\
\end{aligned}$$ 
where, as before, $0\leq a\leq |Q|-1$. A routine check shows that this is well defined and that $\iota(g)k:Y_1\to Y_2$ cyclically preserves  the order between $Y_1$ and $Y_2$, in other words, that $\iota(g)k\in T_r(\Sigma)$. From this we deduce that $\tau: C_{T_r(Q)}(Q) \epi T_l(\Sigma)$ is an epimorphism.

Finally, recall that by the proof of Theorem \ref{centralizer} and Remark \ref{kernel}, $\text{Ker}\tau$ is the union of its subgroups of the form
$$K_Z=\{h\in C_{G_r(\Sigma)}(Q)\mid Zh=Z\text{ and $h$ fixes the $Q$-orbits setwise}\}\cong L^m$$
for each $Q$-expansion $Z$ of $Y$.  Moreover, as $Q$ is transitive and regular, $L=C_{S_{w_1}}(\phi(Q))\cong Q.$ 
The observation before this Theorem implies that the finite group $K_Z\cap T_r(\Sigma)$ must  in fact act in the same way as $Q$ acts on $Z.$ Hence 
$$K_Z\cap T_r(\Sigma)=Q.$$
Note, that under the isomorphism $K_Z\cong Q^m$, this corresponds to the diagonal subgroup of $Q^m$. That the extension is central now follows immediately.
\end{proof}


\begin{remark} For $U_r(\Sigma)=V_{2,1},$ the Higman algebra and $T_r(\Sigma)=T,$ the original Thompson-group $T$ this reproves \cite[Theorem 7.1.5]{matuccithesis}. 
\end{remark}

\section{Finiteness conditions in Bredon cohomology}

\noindent In this section we collect all necessary background on Bredon cohomological finiteness conditions and also prove an analogue to Bieri's criterion for $\FP_n.$ 

\noindent
Let $\mathcal X$ denote a family of  subgroups of a given group $G$.
In Bredon cohomology, the group
$G$ is replaced by the orbit category $\mathcal O_{\mathcal X}G$.  The category $\OXG$ has as objects the transitive $G$-sets with stabilisers in $\mathcal X$. Morphisms in $\OXG$ are $G$-maps between those $G$-sets.  Modules over the orbit category, called $\OXG$-modules are contravariant functors from the orbit category to the category of abelian groups.  Exactness is defined pointwise: a sequence
 $$A\to B\to C$$
  of $\OXG$-modules
is exact at $B$ if and only if
$$A(\Delta)\to B(\Delta)\to C(\Delta)$$
is exact at $B(\Delta)$ for every transitive $G$-set $\Delta$.

\noindent
The category $\OXG$-$\Mod$ of $\OXG$-modules has  enough projectives, which are constructed as follows:
For any $G$-sets $\Delta$ and $\Omega$, denote by $[\Delta,\Omega]$ the set of $G$-maps from $\Delta$ to $\Omega$.
Let $\Z[\Delta,\Omega]$ be the free abelian group on $[\Delta,\Omega]$.
One now obtains an   $\OXG$-module $\Z[\blah,\Omega]$ by fixing $\Omega$ and letting $\Delta$ range over the transitive $G$-sets with stabilisers in $\mathcal X.$ A Yoneda-type argument, see \cite{mislinsurvey}, yields that these modules are free.  In particular, the modules $P^K(\blah)  = \Z[\blah, G/K]$ for $K \in {\mathcal X}$ are free and  can be viewed as the building blocks for free $\OXG$-modules.  Projective modules are now defined analogously to the ordinary case.
The trivial $\OXG$-module, denoted $\Z(\blah)$, is the constant functor $\Z$ from $\OXG$ to the category of abelian groups.

\medskip\noindent Bieri \cite{Bieribook} gives criteria for a $\Z G$-module to be of type $\FP_n$ involving certain $Ext$- and $Tor$-functors  commuting with exact colimits and direct products respectively. In this section we prove that those criteria can also be  used for Bredon cohomology.
The Bredon cohomology functors $\EXT^*_\XX (M,-)$ are defined as derived functors of
$\Hom_\XX(M,-)$.  In particular, let $M(\blah) \in \OXG$-$\Mod$ be a contravariant $\OFG$-module admitting a projective resolution $P_\ast(\blah) \epi M(\blah).$ Then,  for each $N(\blah) \in \OXG$-$\Mod$,
$$ \EXT^*_\XX (M,N) =H_*(mor(P_*,N)).$$

 One can also define Bredon homology functors $\TOR_*^\XX(-,M)$.
 In particular, by analogy with the contravariant case, one can define covariant $\OXG$-modules, or just comodules for short. The category of covariant $\OXG$-modules, denoted $\Mod$-$\OXG$, behaves just as expected.  For example, we have short exact sequences and enough projectives as above. In particular, the building blocks for projective modules in $\Mod$-$\OXG$ are the covariant functors $P_K(\blah)=\Z[G/K, \blah]$ for subgroups $K \in\XX.$ Let $M(\blah) \in \OXG$-$\Mod$ be as before. Then Bredon homology functors are the derived functors of $\blah\otimes_\XX M,$ i.e., for any $L(\blah) \in \Mod$-$\OXG$,
 $$\Tor^\XX_*(L,M) =H_*(L \otimes_\XX P_*).$$

 For detail on these definitions including the categorical tensor product and Yoneda-type isomorphism the reader is referred to \cite{nucinkis04}.  In particular, $\Tor^\XX_*(\blah, M)$ can be calculated using flat resolutions of $M(\blah).$

\medskip\noindent The category of $\OXG$-modules, as an abelian category, has well defined colimits and limits and in particular coproducts and products. We say a functor
$$T: \OXG\mbox{-}\Mod \to {\mathcal A}b$$
commutes with exact colimits, denoted here by $\colim$, if, for every directed system $(M_\lambda)_{\lambda\in \Lambda}$ of $\OXG$-modules, the natural map
$$\colim T(M_\lambda) \to T(\colim M_\lambda)$$
is an isomorphism. Analogously, we say a functor
$$S: \Mod\mbox{-}\OXG \to {\mathcal A}b$$
commutes with exact limits, denoted here by $\invlim$, if, for every inverse system $(N_\lambda)_{\lambda\in \Lambda}$ of $\OXG$-comodules, the natural map
$$S(\invlim N_\lambda) \to  \invlim S(N_\lambda)$$
is an isomorphism.

\noindent We say an $\OXG$-module $M(\blah)$ is finitely generated if there is a finitely generated free module mapping onto it. In particular, there is a $G$-finite $G$-set $\Delta$ such that $\Z[\blah, \Delta] \epi M(\blah)$ (here we are extending the notation $\Z[\blah,\Delta]$ to non transitive sets in the obvious way).

\begin{lemma}\label{fgcolimlemma}
Let $M$ be an $\OXG$-module. Then $M$ is the direct colimit of its finitely generated submodules.
\end{lemma}

\proof This follows from \cite[\S 9.19]{lueckbook}. \qed

\medskip\noindent
The notions of type Bredon-$\FP$, Bredon-$\FP_n$ and Bredon-$\fpinfty$
are defined in terms of projective resolutions  over $\OXG$ analogously to
the classical notions of type $\FP$, $\FP_n$ and $\fpinfty$.

\begin{proposition}\label{ufpnprop}
Let $A$ be an $\OXG$-module of type Bredon-$\FP_n$, $0\leq n\leq \infty.$ Then
\begin{enumerate}
\item For every exact limit, the natural homomorphism
$$\TOR_k^\XX(\invlim N_*, A) \to \invlim\TOR_k^\XX(N_*,A)$$
 is an isomorphism for all $k \leq n-1$ and an epimorphism for $k=n.$
\item For every exact colimit, the natural homomorphism
$$\colim\EXT^k_\XX(A,M_*) \to \EXT^k_\XX(A,\colim M_*)$$
 is an isomorphism for all $k \leq n-1$ and a monomorphism for $k=n.$
\end{enumerate}
\end{proposition}

\proof  The proof goes completely analogously to that of Bieri \cite[Proposition 1.2]{Bieribook}. It relies on the Yoneda isomorphisms, i.e that $N \otimes_\XX \Z[\blah ,G/K] \cong N(G/K)$ and $\Hom_\XX( \Z[\blah, G/K], M) \cong M(G/K)$, the fact that  $\invlim$ and $\Hom_\XX(-, M)$ commute with finite direct sums and that $\invlim$ and $\colim$ are exact and hence commute with the homology functor.    \qed

\medskip\noindent Bieri's argument can be carried through completely for Bredon-Ext and Bredon-Tor functors.
\begin{theorem}\label{ufpnext}
Let $A$ be an $\OXG$-module. Then the following are equivalent:
\begin{enumerate}
\item $A$ is of type Bredon-$FP_n.$
\item For every exact colimit, the natural homomorphism
$$\colim\EXT^k_\XX(A,M_*) \to \EXT^k_\XX(A,\colim M_*)$$
 is an isomorphism for all $k \leq n-1$ and a monomorphism for $k=n.$
\item  For the direct limit of any directed system of $\OXG$-modules $M_*$ with $\colim M_* =0,$ one has $\colim \EXT^k_\XX(A,M_*)=0,$ for all $k \leq n.$
\end{enumerate}
\end{theorem}

\proof  The implications $(i)\implies (ii) \implies (iii)$ are either obvious or follow from Proposition \ref{ufpnprop}. Every $\OXG$-module is the directed colimit of finitely generated submodules, Lemma \ref{fgcolimlemma}, and hence $(iii)\implies (i)$ is proved completely analogously to \cite[Theorem 1.3 $(iiib) \implies (i)$]{Bieribook}.  \qed

\medskip

\begin{theorem}\label{ufpntor} Let $A$ be an $\OXG$-module. Then the following are equivalent:
\begin{enumerate}
\item $A$ is of type Bredon-$\FP_n.$
\item For every exact limit, the natural homomorphism
$$\TOR_k^\XX( \invlim N_*, A) \to \invlim\TOR_k^\XX(N_*, A)$$
 is an isomorphism for all $k \leq n-1$ and an epimorphism for $k=n.$
\item For any $K\in\XX$ consider any arbitrary  direct product $\prod_{\Lambda_K} \Z[G/K, \blah]$. Then the natural map
$$\TOR_k^\XX(\prod_{K\in\XX} \prod_{\Lambda_K} \Z[G/K, \blah], A) \to \prod_{K\in\XX} \prod_{\Lambda_K}\TOR_k^\XX(\Z[G/K, \blah], A)$$
 is an isomorphism for all $k \leq n-1$ and an epimorphism for $k=n.$

\end{enumerate}
\end{theorem}

\proof The implications $(i)\implies (ii) \implies (iii)$  are again either obvious or consequence of Proposition \ref{ufpnprop}.

\noindent $(iii) \implies (i): $ The proof  is in the same spirit as Bieri's proof. We begin by letting $n=0$ and claim that $A$ is finitely generated as an $\OXG$-module. As an index set we take $\prod_{K\in\XX} A(G/K)$ and consider $\prod_{K\in\XX} \prod_{a \in A(G/K)} \Z[G/K ,\blah].$ By $(iii),$ the natural map
$$ \mu: (\prod_{K\in\XX} \prod_{A(G/K)} \Z[G/K, -]) \otimes_\XX A(-) \to \prod_{A(G/K)}A(G/K) $$
is an epimorphism.	Let $c$ be the element with $\mu(c)=\prod_{K\in\XX} \prod_{a\in A(G/K)}a $. Then $c$  is of the form
$$c = \sum_{i=1}^l (\prod_{K\in\XX} \prod _{A(G/K)}f_i^{a,K})\otimes b_i,$$
for certain subgroups $H_1,\ldots,H_l\in\XX$ and elements $b_i\in A(G/H_i)$. Here, $f_i^{a,K} \in\Z[G/K,G/H_i]$.
Now we claim that there is an epimorphism
$$\tau:\bigoplus_ {i=1}^l\Z[\blah,G/H_i]\twoheadrightarrow A$$
given by $\tau(f):=f^*(b_i)\in A(G/K)$ whenever $f\in\Z[G/K,G/H_i].$ Observe that this is well defined. In particular it is functorial.
To prove the claim, take any $K\in\XX$ and any $a\in A(G/K)$. Note that
$$\mu(c)=\sum_{i=1}^l\prod_{K\in\XX} \prod_{a\in A(G/K)} (f_i^{a,K})^*(b_i)=\prod_{K\in\XX} \prod_{a\in A(G/K)}\sum_{i=1}^l(f_i^{a,K})^*(b_i)$$
so the fact that $c$ maps onto the diagonal means that
$$a=\sum_{i=1}^l(f_i^{a,K})^*(b_i)=\tau(\sum_{i=1}^lf_i^{a,K}).$$

\noindent The case $n \geq 1$ is now done analogously to \cite[Theorem 1.3]{Bieribook} using a diagram chase. \qed

\begin{remark}  For $n\geq 1$, condition (iii) is equivalent to the following, which in ordinary homology is often referred to as the Bieri-Eckmann criterion for $\FP_n:$   For every  subgroup $K \in \XX$  consider an arbitrary  direct product $\prod_{ \Lambda_K} \Z[G/K, \blah]$. Then the natural map
 $$ (\prod_{K\in\XX}  \prod_{ \Lambda_K} \Z[G/K, \blah ]) \otimes_\XX A(\blah) \to \prod_{K\in\XX} \prod_{ \Lambda_K} A(G/K)$$ is an isomorphism and $\Tor^\XX_k(\prod_{K\in\XX} \prod_{ \Lambda_K} \Z[G/K, \blah ], A) =0,$ for all $1 \leq k \leq n-1.$ We call this condition the global Bieri-Eckmann criterion for Bredon homology.

 We say a group satisfies the local Bieri-Eckmann criterion for Bredon cohomology if, for any $K$   and direct product as before, the natural map
 $$ (\prod_{ \Lambda_K} \Z[G/K, \blah ] )\otimes_\XX A(\blah) \to  \prod_{ \Lambda_K} A(G/K)$$ is an isomorphism and $\Tor^\XX_k( \prod_{ \Lambda_K} \Z[G/K, \blah ], A) =0$ for all $1 \leq k \leq n-1.$
 \end{remark}

\section{Classifying spaces with finite isotropy}

\noindent In this section we shall restrict ourselves to the family ${\mathcal F}$ of all the finite subgroups of $G$.

To stay in line with notation previously used, we say a module is of type $\UFP_\infty$  if it is of type Bredon-$\fpinfty$ with respect to  $\mathcal F$. The notions of $\UFP_n$ and $\UFP$ are defined analogously.
For Bredon cohomology with respect to $\mathcal F$ there is a good algebraic description for modules of type $\UFP_n$. For the original approach via classifying spaces, see
\cite{lueck}.

\begin{theorem}\cite{kmn}\label{ufpn}
Let $G$ be a group having finitely many conjugacy classes of finite subgroups. Then
an $\OFG$-module $M(\blah)$   is of type $\UFP_n$ if and only if  $M(G/K)$ is of type $\FP_n$ as a $\Z(WK)$-module for each finite subgroup $K$ of $G$.
\end{theorem}

It was also shown, \cite{kmn}, that a group $G$ is of type $\UFP_0$ if and only if $G$ has finitely many conjugacy classes of finite subgroups. Hence we have the following corollary:

\begin{corollary}\cite{kmn}\label{Gufpn}
A group $G$ is of type $\UFP_n$ if and only if $G$ has finitely many conjugacy classes of finite subgroups and $C_G(K)$ is of type $\FP_n$ for every finite subgroup $K$ of $G$.
\end{corollary}

\noindent Recall that we say a group $G$ is of Bredon-type $\FP_n$ if the trivial module $\Z(\blah)$ is of type $\FP_n$ as an $\OXG$-module.  We can, of course rephrase Theorems \ref{ufpnext} and  \ref{ufpntor} in terms of Bredon-cohomology and Bredon-homology replacing the module $A(\blah)$ with $\Z(\blah)$, $\EXT^*_\XX(A,-)$ with $\Ho^*_\XX(G,-)$ and $\Tor^\XX_*(-,A)$ with $\Ho^\XX_*(G,-).$


\medskip\noindent We shall now weaken the hypothesis on the conjugacy classes of finite subgroups:

\begin{definition}
We say a group is of type quasi-$\UFP_n$ if, for each finite subgroup $K$ of $G$ there are finitely many conjugacy classes of subgroups isomorphic to $K$ and the Weyl-groups $WK:=N_G(K)/K$ are of type $\FP_n.$
\end{definition}

\noindent
Note that a group of type quasi-$\UFP_n$ with a bound on the orders of the finite subgroups is of type $\UFP_n$.

\noindent
Let $k$ be a positive integer.  We denote by $\Z_k(\blah)$ the $\OFG$-module defined by
$$\Z_k(G/H) = \begin{cases} \Z \mbox{ if } |H|\leq k \\
                                                    0 \mbox{ otherwise, } \end{cases}$$
                                                    together with the obvious morphisms.

\begin{lemma}\label{quasifp0}
A group $G$ is of type quasi-$\UFP_0$ if and only if, for each $k \geq 1$, the module $\Z_k(\blah)$ is finitely generated. Moreover, in that case, the finite $G$-set $\Delta_k$ with  $\Z[ -, \Delta_k] \epi \Z_k(\blah)$ can be chosen to have stabilisers of order bounded by $k$.
\end{lemma}

\proof Suppose $G$ is of type quasi-$\UFP_0.$ Take
$$\Delta_k = \bigsqcup_{|H|\leq k, { \mbox{up to $G$-conj.}}} G/H.$$
 This is a $G$-finite $G$-set with stabilisers of order bounded by $k$ and
$\Z[ -, \Delta_k] \epi \Z_k(\blah).$  \noindent
For the converse, we need to show that, for each finite subgroup $K$, there are only finitely many conjugacy classes of subgroups of order bounded by $k=|K|$. Let $\Delta_k$ be the finite $G$-set with $\Z[ -, \Delta_k] \epi \Z_k(\blah)$ and take any finite subgroup $H$ of $G$ with $|H| \leq k$. Hence $\Z_{k}(G/H)\cong \Z \neq 0.$ Since the map $\Z[G/H,\Delta] \epi \Z_{k}(G/H)$ is onto, it follows that $\Z[G/H,\Delta] \neq 0$ and hence $H$ has to be subconjugated to one of the finitely many stabilisers of $\Delta.$     \qed

\medskip\noindent Note that finitely generated $\OFG$-modules are precisely those of type $\UFP_0.$ Fix an integer $k \geq 1$ and
let $M_k(-)$ be an $\OFG$-module such that $M_k(G/L)=0$ whenever $|L|>k.$ Suppose $M_k(-)$ is finitely generated. Then there exists a $G$-finite $G$-set $\Delta$ with stabilisers of order $\leq k$ and a short exact sequence of $\OFG$-modules
$$N_k(-) \mono \Z[-,\Delta] \epi M_k(-)$$
with the property that $N_k(G/L)=0$ for all finite subgroups $L$ with $|L|>k.$

\begin{proposition}\label{quasifpn}
A group $G$ is of type quasi-$\UFP_n$ if and only if, for each integer $k \geq1,$ the $\OFG$-module $\Z_k(\blah)$ is of type $\UFP_n.$
\end{proposition}

\proof  The "if"-direction follows from Lemma \ref{quasifp0}, Theorem \ref{ufpn} and the definition as $\Z_{|K|}(G/H)$ is of type $\FP_n$ as a $WH$-module for each $|H| \leq |K|.$

\noindent
Now suppose $G$ is of type quasi-$\UFP_n.$ For each $k \geq 1$ we construct a projective resolution of $\Z_k(\blah)$ which is finitely generated in dimensions up to $n$; note that we may assume $n>0$.
By Lemma \ref{quasifp0} and the above remark we have a short exact sequence
$$C_0(\blah) \mono \Z[-,\Delta_0] \epi \Z_k(\blah)$$
with $\Delta_0$ a $G$-finite $G$-set and $C_0(G/L)=0$ for all $|L|>k.$
We claim that $C_0(-)$ is a finitely generated $\OFG$-module.

 We know that there are finitely many conjugacy classes of subgroups of order bounded by $k$. Let $H$ be one of those. As $\Delta_0$ is $G$-finite, the $WH$-module $\Z[G/H,\Delta_0]$ is of type $\FP_\infty$. This is a consequence of the fact that for any $K$,
 $\Z[G/H,G/K]$ is a sum of exactly $|\{x\in N_G(H)\backslash G/K\mid H^{x^{-1}}\leq K\}|$ $WH$-modules, which are of type $\FP_\infty$. As $K$ is finite, this sum must also be finite.
So evaluating  the previous short exact sequence at $G/H$, we see that the $WH$-module $C_0(G/H)$ is of type $\FP_{n-1}$ and in particular, finitely generated. Fix a finite $WH$-generating set $X_H$ for $C_0(G/H)$. Then the $\OFG$-set formed by the union of all those $X_H$ where $H\in\text{Stab}\Delta_0,$ generates $C_0$.


We can now proceed to construct the desired resolution by using the remark before Proposition \ref{quasifpn}. \qed

\begin{theorem}\label{bierieckmann}
Let $G$ be of type quasi-$\UFP_n,$  where $n \geq 1.$ Then $G$ satisfies the local Bieri-Eckmann criterion for Bredon homology.

 \end{theorem}

\proof
It follows from the definition of the modules $\Z_k(\blah)$ that
$$\Z(-) = \colim_{k \in \N} \Z_k(-).$$
In the category of $\OFG$-modules the construction of a free module mapping onto a given one is functorial. Hence, we can get a direct colimit of free resolutions $\colim_{k \in \N}(F_{*,k}(-) \epi \Z_k(-))=F_*(-) \epi \Z(-)$, which gives us a flat resolution of $\Z(-).$ For details the reader is referred to \cite[Lemma 3.4]{nucinkis04}.
Hence
\begin{eqnarray*}
\Ho^{\mathcal F}_k(G,  \prod_{ \Lambda} \Z[G/K, \blah ])  = \Ho_*(  \prod_{ \Lambda} \Z[G/K, \blah ]\otimes_{\mathcal F} F_*(-))& \\
= \Ho_k( \prod_{ \Lambda} \Z[G/K, \blah ]\otimes_{\mathcal F}\colim_{k \in \N}F_{*,k}(-)) &\\
= \colim_{k \in \N}\Ho_k(  \prod_{ \Lambda} \Z[G/K, \blah ]\otimes_{\mathcal F} F_{*,k}(-))& \\
=  \colim_{k \in \N} \Tor_k(\prod_{ \Lambda} \Z[G/K, \blah ],\Z_k(-)) = 0&,\end{eqnarray*}
where the last line follows from Proposition \ref{quasifpn} and Theorem \ref{ufpntor}.
The first assertion follows by a similar argument.

\qed

\medskip\noindent For each $k \geq 1$ we consider the family ${\mathcal F}_k$ and the orbit category $\OFKG$.
For a given positive integer  $k$ the family ${\mathcal F}_k$ consists of all subgroups $H$ of $G$ with $|H| \leq k.$ By using the arguments of the proofs of Lemma \ref{quasifp0} and Proposition \ref{quasifpn} we can show:

\begin{proposition}\label{quasifk} A  group is of type quasi-$\UFP_n$ if and only if it is of type Bredon-$\FP_n$ over $\OFKG$ for each $k$.
\end{proposition}


\noindent We can also rephrase Theorems and \ref{ufpnext} and \ref{ufpntor}:

\begin{corollary}\label{quasihom}
Let $G$ be a group. Then the following are equivalent:
\begin{enumerate}
\item $G$ is of type quasi-$\UFP_n.$
\item For every exact colimit and any $k,$ the natural homomorphism
$$\colim\Ho^l_\FK(G,M_*) \to \Ho^l_\FK(G,\colim M_*)$$
 is an isomorphism for all $l \leq n-1,$ and a monomorphism for $l=n.$
\item For any $k$ and any $K\in\FK$ consider an arbitrary  direct product $\prod_{\Lambda_K} \Z[G/K, \blah]$. Then the natural map
$$\Ho_l^\FK(\prod_{K\in\FK} \prod_{\Lambda_K} \Z[G/K, \blah], A) \to \prod_{K\in\FK} \prod_{\Lambda_K}\Ho_l^\FK(\Z[G/K, \blah], A)$$
 is an isomorphism for all $l \leq n-1$ and an epimorphism for $l=n.$
\end{enumerate}
\noindent One may also add the statements analogous to \ref{ufpnext} ii) and \ref{ufpntor} ii). Note also that  for $n\geq 1$ the above is equivalent to:
\begin{itemize}
\item[(iv)] For any $k$, any $K\in\FK$ and any arbitrary  direct product $\prod_{ \Lambda_K} \Z[G/K, \blah],$ the natural map

 $$ \Z_k(\blah)  \otimes_\FK \prod_{K\in\FK}\prod_{\Lambda_K} \Z[G/H, \blah ] \to \prod_{K\in\FK}\prod_{\Lambda_K} \Z_k$$ is an isomorphism and $\Ho^\FK_l(G, \prod_{K\in\FK}\prod_{\Lambda_K} \Z[G/H, \blah ]) =0$ for all $1 \leq l \leq n-1.$

\end{itemize}
\end{corollary}

\begin{definition}
We say a group $G$ is of type quasi-$\UF_\infty$ if for all positive integers $k$, $G$ admits a finite type model for $E_\FK G.$
\end{definition}

\noindent Analogously to the algebraic case, any group of type quasi-$\UF_\infty$, which has a bound on the orders of the finite subgroups, is of type $\UF_\infty.$

\noindent L\"uck's Theorem \cite[Theorem 4.2]{lueck} (Theorem \ref{lueckstheorem}) goes through for arbitrary families of finite subgroups. Hence combining Theorems \ref{lueckmeintrupp} and \ref{lueckstheorem} yields:

\begin{proposition}\label{uf}
A group $G$ is of type quasi-$\UF_\infty$ if and only if $G$ is of type quasi-$\UFP_\infty$ and $G$ and all centralisers $C_G(K)$ of finite subgroups are finitely presented. \qed
\end{proposition}

\noindent We can now prove what is largely equivalent to Proposition \ref{quasifk}:

\begin{theorem}
A group $G$ is of type quasi-$\UF_\infty$ if and only if it admits a model for $\eg$, which is the mapping telescope of finite type models for $E_\FK G$ for each $K \in {\mathcal F}.$
\end{theorem}

\proof The "if'-direction follows directly from the definition. Now suppose we have finite type models $X_K$ for $E_\FK G$ for all $K \in {\mathcal F}.$ For each $H \leq K$ the universal property for classifying spaces for a family yields $G$-maps $\nu_H^K: X_H \to X_K.$ Now the mapping telescope yields a $G$-CW-complex X, for which $X^K$ is contractible for all $K \in {\mathcal F}$ and empty otherwise.  \qed

\section{Bredon cohomological finiteness properties for generalised Thompson-Higman groups}

\noindent We can now prove

\begin{theorem}\label{Tmain2} Let $U_r(\Sigma)$ be a Cantor-Algebra with order preserving $\Sigma$. Then the following conditions are equivalent for $1\leq r\leq d$:
\begin{itemize}

\item[i)] $T_r(\Sigma)$ is quasi-$\underline{\FP}_\infty$.

\item[ii)] $T_l(\Sigma)$ is of (ordinary) type  $\FP_\infty$ for any $1\leq l\leq d$ such that $\text{gcd}(l,d)\mid r$.

\end{itemize}
\end{theorem}
\begin{proof}  Assume that i) holds and take $1\leq l\leq d$ with $\text{gcd}(l,d)\mid r$. This condition implies that there is some positive integer $w$ with $lw\equiv r$ mod $d$. Then we may choose an admissible subset $A\subseteq U_r(\Sigma)$ of cardinality $lw$ and consider the subgroup $Q$ of $T_r(\Sigma)$ defined by cyclic permutations of $A$ on $l$ orbits all of length $w$. By Theorem \ref{Textension}, $C_{T_r(\Sigma)}(Q)$ is an extension of a finite group by $T_l(\Sigma)$ so this last group must be of type $\FP_\infty$.

Now assume ii). Observe,  that for any finite subgroup $Q$ of cardinality $w$ and any admissible subset $Y_1$ fixed by $Q$, the observation before Theorem \ref{Textension} implies that for certain $l$, $|Y_1|=lw\equiv r$ mod $d$ thus $\text{gcd}(l,d)\mid r$. This together with Theorem \ref{Textension} implies that $T_r(\Sigma)$ is quasi-$\underline{\FP}_\infty$.

\end{proof}

\noindent We also have the same result for finiteness conditions on classifying spaces

\begin{theorem}\label{Tmain} Let $U_r(\Sigma)$ be a Cantor-Algebra with order preserving $\Sigma$. Then the following conditions are equivalent:
\begin{itemize}

\item[i)] $T_r(\Sigma)$ is quasi-$\underline{\FFF}_\infty$.

\item[ii)] $T_l(\Sigma)$ is of (ordinary) type  $\FFF_\infty$ for any $1\leq l\leq d$ such that  $\text{gcd}(l,d)\mid r$.
\qed

\end{itemize}

\end{theorem}

\begin{proof} This follows from Theorem \ref{Textension} exactly as Theorem \ref{Tmain2}.
\end{proof}

\begin{corollary}
Let $U_r(\Sigma)$ be a Higman algebra. Then 
$T_{n,r}=T_r(\Sigma)$ is quasi-$\underline{\FFF}_\infty$.
\end{corollary}

\begin{proof}
This follows directly from \cite{brown} and Theorem 
\ref{Tmain}.
\end{proof}

\begin{corollary} Let $U_r(\Sigma)$ be a Brown-Stein algebra. Then 
$T=T_r(\Sigma)$ is quasi-$\underline{\FFF}_\infty$.
\end{corollary}

\begin{proof} This is a consequence of Theorem 
\ref{Tmain} and  \cite[Theorem 2.5]{stein} where it is proven that $F_r(\Sigma)$ is finitely presented and of type $\FP_\infty$ for any $r$. Stein's argument carries over to $G$ and $T$,  \cite{stein}.
\end{proof}

\begin{conjecture} Let $U_r(\Sigma)$ be a Cantor-algebra. Then
\begin{enumerate}
\item $G_r(\Sigma)$ is quasi-$\underline{\FP}_\infty$ if and only if $G_l(\Sigma)$ is of (ordinary) type ${\FP}_\infty$ for any $1\leq l\leq d.$ 
\item $G_r(\Sigma)$ is quasi-$\underline{\FFF}_\infty$ if and only if $G_l(\Sigma)$ is of (ordinary) type ${\FFF}_\infty$ for any $1\leq l\leq d.$ 
\end{enumerate}
\end{conjecture}

\begin{remark} Our description of centralisers of finite subgroups implies the \lq\lq only if" part of this conjecture. To see this, assume that $G_r(\Sigma)$ is of type quasi-$\underline{\FP}_\infty$ and choose for any $1\leq l\leq d$ positive integers $n,s$ with $n\geq 3$ such that $ln+s\equiv r$ mod $d$. Then there is some admissible subset $Y\in\frak A_r(\Sigma)$ of cardinality precisely $ln+s$ and we may consider the finite group $Q\cong S_n\leq G_r(\Sigma)$ defined by the action on $l$ orbits of $n$ elements as the natural representation
of $S_n$ and acting trivially on the remaining $s$ elements of $Y$. Then Theorem \ref{centralizer} implies, using the same notation here, that 
$$C_{G_r(\Sigma)}(Q)\cong H_l\times H_s.$$
Hence both $H_l$ and $H_s$ are of type $\FP_\infty$. Moreover, $H_l$ is an extension
$$K_1\rightarrowtail H_l\twoheadrightarrow G_l(\Sigma),$$
where by Remark \ref{kernel} $K_1$ is a direct limit of products of $L=C_{S_n}(S_n)=1$. Thus, in this case, $L=K_1=1$ implying that $G_l(\Sigma)$ is of type $\FP_\infty$.

\end{remark}

\begin{remark}\label{remark1} By \cite{brown}, Proposition 4.1, $F_{n,r}\cong F_{n,s}$, for any $r,s$. However, this is false for the groups $G$, in fact $G_{n,r}\cong G_{n,s}$ implies $\text{gcd}(n-1,r)=\text{gcd}(n-1,s)$ (\cite{higman} Theorem 6.4). Recently, Pardo founded \cite{pardo} that the converse also holds true (see also \cite{dicksconcha}).
\end{remark}

\end{document}

\noindent Recall that a group is of type $\UFP_n$ if the trivial module $\uz$ is of type $\UFP_n$ as an $\OFG$-module. It was shown in \cite{kmn} that a group is of type $\UFP_0$ if and only if it has finitely many conjugacy classes of finite subgroups. Hence $\uz$ is a finitely generated $\OFG$-module if and only if $G$ has finitely many conjugacy classes of finite subgroups and we have the following corollary to Theorem \ref{ufpnext}:

\begin{corollary}\label{ufpncoho}
For any group $G$ the following are equivalent:
\begin{enumerate}
\item $G$ is of type $\UFP_n.$
\item  For every exact colimit, the natual homomorphism
$$\colim\Ho^k_{\mathcal F}(G,M_*) \to \Ho^k_{\mathcal F}(G,\colim M_*)$$
 is an isomorphism for all $k \leq n-1$ and a monomorphism for $k=n.$
\item  For the direct limit of any directed system of $\OFG$-modules $M_*$ with $\colim M_* =0$ one has $\colim \Ho^k_{\mathcal F}(G,M_*)=0$ for all $k \leq n.$
\end{enumerate}
\end{corollary}

\noindent Also, a group is of type $\UFP_n$ if and only if it has finitely many conjugacy classes of finite subgroups and the Weyl-group $WK$ is of type $\FP_n$ for each finite subgroup $K$ of $G$ \cite{kmn}.
Hence we have the following corollary to Theorem \ref{ufpntor}:

\begin{corollary}\label{ufpnho}
Let $G$ be a group with finitely many conjugacy classes of finite subgroups. Then the following are equivalent:
\begin{enumerate}
\item $G$ is of type $\UFP_n.$
\item For every exact limit, the natural homomorphism
$$\Ho_k^{\mathcal F}(G, \invlim N_*) \to \invlim\Ho_k^{\mathcal F}(G,N_*)$$
 is an isomorphism for all $k \leq n-1$ and an epimorphism for $k=n.$
\item  For every finite subgroup $K$ of $G$ and a arbitrary  direct product $\prod \Z[G/K, \blah]$ the natural map
$$\Ho_k^{\mathcal F}(G, \prod \Z[G/K, \blah]) \to \prod\Ho_k^{\mathcal F}(G,\Z[G/K, \blah])$$
 is an isomorphism for all $k \leq n-1$ and an epimorphism for $k=n.$
\end{enumerate}
\end{corollary}